\documentclass[11pt,a4paper]{article}
\usepackage[top=3cm, bottom=3cm, left=3cm, right=3cm]{geometry}
\usepackage[latin1]{inputenc}
\usepackage{amsmath}
\usepackage{amsthm}
\usepackage{amsfonts}
\usepackage{amssymb}
\usepackage{graphics}
\usepackage{tabu}
\usepackage{breqn}
\usepackage[hang,flushmargin]{footmisc} 

\usepackage{amsmath,amsthm,amssymb,graphicx, multicol, array}
\usepackage{enumerate}
\usepackage{enumitem}
\usepackage{hyperref}
\usepackage{setspace}
\usepackage{xcolor}
\usepackage{currfile}

\usepackage{tabto}

\usepackage{pgfplots}
\pgfplotsset{width=10cm,compat=1.9}

\newtheorem{thm}{Theorem}[section]

\newcommand{\C}{\mathbb{C}}

\usepackage{fancyhdr}
\title{Simple Zeros Of The Zeta Function}
\date{}
\author{N. A. Carella}

\begin{document}
\thispagestyle{empty}
\date{}

\maketitle

\textbf{\textit{Abstract}:} This note studies the Laurent series of the inverse zeta function $1/\zeta(s)$ at any fixed nontrivial zero $\rho$ of the zeta function $\zeta(s)$, and its connection to the simplicity of the nontrivial zeros.  \let\thefootnote\relax\footnote{ \today \date{} 
}

\section{Introduction} \label{s111}
The set of zeros of the zeta function $\zeta(s)=\sum_{n\geq 1}n^{-s}, \Re e(s)>1$ is defined by $\mathcal{Z} = \# \{ s \in \mathbb{C} : \zeta(s) = 0 \}$.  This subset of complex numbers has a disjoint partition as $\mathcal{Z}=\mathcal{Z}_T \cup \mathcal{Z}_N$, the subset of trivial zeros and the subset of nontrivial zeros. \\

The subset of trivial zeros of the zeta function is completely determined; it is the subset of negative even integers $ \mathcal{Z}_T=\{-2n:\zeta(s)=0 \text{ and } n \geq1\}$ . 
The trivial zeros are extracted using the symmetric functional equation
\begin{equation}
\pi^{-s/2} \Gamma(s/2) \zeta(s)=\pi^{-(1-s)/2} \Gamma((1-s)/2) \zeta(1-s) ,
\end{equation}
where $\Gamma(t)=\int_{0}^{\infty}x^{t-1}e^{-x}d x$ is the gamma function, and $s\in \mathbb{C}$ is a complex number. Various levels of explanations of the functional equation are given in \cite[p. 222]{NW00}, 
\cite[p. 328]{MV07}, \cite[p. 13]{IA03}, \cite[p. 55]{KA93}, \cite[p. 8]{IA03}, \cite{TE86}, \cite{LMFDB}, and other sources. In contrast, the subset of nontrivial zeros 
\begin{equation}
\mathcal{Z}_N=\{s \in \mathbb{C}:\zeta(s)=0 \text{ and } 0 < \Re e(s) <1\}
\end{equation}
and myriad of questions about its properties remain unknown. The nontrivial zeros are computed using the Riemann-Siegel formula, see \cite[Eq. 25.10.3]{DLMF}. Some recent development 
in this area appears in \cite{HG11} and \cite{HG16}. \\
		
Let $N(T)=\#\{ \rho=\sigma +it:\zeta(\rho)=0, \quad 0<\sigma<1, \text{ and }|t| \leq T \}$ be the counting function of the subset of nontrivial zeros, and let
\begin{equation}
N_0(T)=\#\{ \rho=\sigma +it:\zeta(\rho)=0, \text{ and }|t| \leq T \}
\end{equation}
and\\ 
\begin{equation}
N_1(T)=\#\{ \rho=\sigma +it:\zeta(\rho)=0, \quad \zeta^{''}(\rho)\ne0,\text{ and }|t| \leq T \}
\end{equation}
be the counting functions of the subsets of nontrivial distinct zeros, and simple zeros respectively. It is clear that $N_1(T) \leq N_0(T)\leq N(T)$, and the RH predicts that $ N_1(T)= N_0(T)= N(T)$, a proof conditional on the RH appears in \cite[p. 5]{OT85}, and \cite{TE86}. The \textit{Riemann Hypothesis} (RH) claims that the nontrivial zeros of the zeta function $\zeta(s)$ are of the form $\rho, t \in \mathbb{R}$. \\
		
A very recent result in \cite{BH13} shows, that assuming the RH, the subset of nontrivial simple zeros has positive density in the set of nontrivial zeros. Specifically,
\begin{equation}
N_1(T) \geq (19/29)N(T)   \qquad \text{  and  }  \qquad    N_0(T) \geq .84665N(T) 
\end{equation}
for all large real number $T\geq 1$. The earliest results on the number of the nontrivial simple zeros are the Selberg result of positive density for nontrivial zeros of odd multiplicities, and the Montgomerry result $N_1(T) \geq (2/3+o(1))N(T) $. The former is based on the Fourier analysis of the convolution of certain functions, refer to \cite{MH72}. A simpler and expanded version of this analysis is given in \cite{GO00}. Other methods and further progress on this problem appear in \cite{CG98}, and the references within. A survey on the theory of the zeros of the zeta function and recent developments appears in \cite{SN10}, and certain associated problems such as the moment of the derivatives are considered in \cite{MN06}. Another related result, given in \cite{SJ02}, claims that the short interval $[T, T + T^{.552}]$ contains a positive proportion of simple zeros and on 
the critical line. This is a vast improvement over the naive estimate $N_1(T+T^.552)-N_1(T) \geq cT^{.552}\log T$, with $c > 0$ constant. \\
			
The objective of this short note is to offer a simpler and indirect approach to the verification of the simplicity of the nontrivial zeros of the zeta function. The main result is the following.			
\begin{thm} \label{thm111.01}  Every zero $\rho = \sigma + it$ of the zeta function $\zeta(s)$ is a simple zero. \end{thm}
			
The proof of this result is an immediate consequence of Theorem \ref{thm222.02} in Section \ref{s222}, which constructs the Laurent series of the inverse zeta function $1/\zeta(s)$. A second independent proof also appears in Theorem \ref{thm666.02} in Section \ref{s666}. A corollary of these results is that the \textit{k}th derivative $\zeta^{(k)}(\rho)\ne0$ for any nontrivial zero $\rho = \sigma + it$, and $k \geq 1$. A couple of related theorems in \cite{LM74} and \cite{YC00} claim that the RH implies that the derivatives $\zeta^{'}(s)\ne0$ and  $\zeta^{''}(s)\ne0$ for any $s = \sigma + it, 0 < \sigma < 1/2$. In general, the zeros of the $k$th derivatives $\zeta^{(k)}(s)$ of the zeta function have complex patterns, confer \cite{BS02} for recent works on this topic. 

\vskip .75 in

\section{The Laurent Series Expansion} \label{s222}
The development of the Laurent series \\
\begin{equation}
\frac{1}{\zeta(s)}=\frac{c_{-m}}{(s-\rho)^m}+ \frac{c_{-m+1}}{(s-\rho)^{m+1}}+\cdots +\frac{c_{-1}}{s-\rho}+c_0+c_1(s-\rho)+c_2(s-\rho)^2+\cdots
\end{equation}
for the inverse zeta function $1/\zeta(s)$ at a fixed nontrivial zero $\rho = \sigma + it$ of the zeta function of multiplicity $m(\rho) = m \geq 1$ is based on an extension of a 
technique used in \cite[p. 206]{CO07} to derive the Laurent series
\begin{equation} \label{el06}
\zeta(s)=\frac{1}{s-1}+\sum_{n\geq 0} \frac{(-1)^n \gamma_n}{n!}(s-1)^n
\end{equation}
of the zeta function at $s = 1$, and the standard practices of the principle of analytic continuation. The first coefficient $\gamma=\gamma_0$ coincides with Euler constant. A 
related analysis yields the Laurent series of the zeta function of a quadratic field $K =\mathbb{Q}(\sqrt{d})$, which has the form:\\
\begin{equation}
\zeta_K(s)=\sum_{(m,n)\ne(0,0)}\frac{1}{Q(m,n)^s}=\frac{2 \pi}{s-1}+\sum_{n\geq 0} \frac{(-1)^n \beta_n}{n!}(s-1)^n,
\end{equation}
where the quadratic form is $Q(x,y)=ax^2+bxy+cy^2$ with $a,c>0$, and $b^2-4ac=-1$. The first coefficient $\beta_0=4 \pi (\gamma+\log (c)/2+\log|\eta(w)|^2)$ is a constant, where 
$\eta(z)=e^{i\pi z/12}\prod_{n\geq 1}\left (1-e^{i \pi nz} \right )$ is the Dedekind eta function, and $w=(b+i)/(2a)$, see \cite[p. 157]{ZD75}. \\

The technique employed here is somewhat related to the Stieltjes-Hermite method, see \cite[p. 160]{NW00}, and \cite[p. 15]{IA12}. Other information on the Stieltjes-Hermite method appear in \cite{AJ11}, \cite{AP85}, \cite{BC88}, \cite{CO07}, \cite{KC11}, and similar sources. As in the case of the series (\ref{el06}), this technique for computing the power series of holomorphic and meromorphic functions, seems to show that any nontrivial zero must have multiplicity $m = 1$ unconditionally. Equivalently, it shows that the nontrivial zeros are simple zeros, unconditionally.

\begin{thm} \label{thm222.02} For each fixed nontrivial zero $\rho = \sigma + it$ of the zeta function $\zeta(s)$, the Laurent series of the inverse zeta function $1/\zeta(s)$ at $\rho = \sigma + it$ has the form
\begin{equation}
\frac{1}{\zeta(s)}=\frac{c_{-1}}{s-\rho}+\sum_{n\geq 0} \frac{(-1)^n \phi_n}{n!}(s-\rho)^n,
\end{equation}
where $c_{-1}$ is the residue at $\rho = \sigma + it$, and the $n$th coefficient $\phi_n$ is given by
\begin{equation}
\phi_n=\sum_{n\geq 1} \left (  \frac{\mu(k) \log( k)^n}{k^{\rho}}- \frac{\log( k+1)^{n+1}-\log( k)^{n+1}}{n+1} \right )
\end{equation}
for $n\geq 0$. This is an analytic function on the domain $D(\rho)=\{ s:0<|s-\rho|<r\}$ of some radius $r > 0$.
\end{thm}
		
\begin{proof} Let $\rho=\sigma_0+it_0$ be a fixed nontrivial zero, $1/2 \leq\sigma_0 < 1$, let $k \geq 1$ be an integer, and consider the sequence of complex valued functions
\begin{equation}\label{el10}
v_k(s)=\frac{\mu(k)}{k^s}-c_{-1}\int_{k}^{k+1} \frac{1}{x^{s-\rho+1}}d x,
\end{equation}
for $s = \sigma + it \in \mathbb{C}$, and the corresponding inequality
\begin{equation}
|v_k(s)|= \left | \frac{\mu(k)}{k^s}-c_{-1}\int_{k}^{k+1} \frac{1}{x^{s-\rho+1}}d x \right | \leq \frac{c}{k^\sigma},
\end{equation}
where $c > 0$ is a constant, and $\Re e(s) = \sigma$.\\

In terms of exponential functions, the $k$th function $v_k(s)$ has the form
\begin{eqnarray}\label{el12}
v_k(s)&=&\frac{\mu(k)}{k^s}+\frac{c_{-1}}{s-\rho}\left (\frac{k}{(k+1)^{s-\rho}}- \frac{1}{k^{s-\rho}} \right )  \nonumber\\
&=&\frac{\mu(\rho)}{k^s}e^{-(s-\rho) \log (k)}+\frac{c_{-1}}{s-\rho}\left ( e^{-(s-\rho) \log (k+1)}-e^{-(s-\rho) \log (k)} \right )\\
&=&\frac{\mu(k)}{k^\rho}  \sum_{n\geq 0}\frac{(-1)^n(s-\rho)^n}{n!} \log (k)^n+        \frac{c_{-1}}{s-\rho}\sum_{n\geq 0}\frac{(-1)^n(s-\rho)^n}{n!}\left ( \log (k+1)^{n}-\log (k)^{n} \right )  \nonumber\\
&=&\frac{\mu(k)}{k^\rho}  \sum_{n\geq 0}\frac{(-1)^n(s-\rho)^n}{n!} \log (k)^n  \nonumber\\
& & \hskip 1 in +     \quad   c_{-1}\sum_{n\geq 0}\frac{(-1)^{n+1}(s-\rho)^n}{(n+1)!}\left ( \log (k+1)^{n+1}-\log (k)^{n+1} \right )  \nonumber.
\end{eqnarray}
The summation index of the second power series was shifted because the first term
\begin{equation}
\frac{(-1)^0 (s-\rho)^0}{0!} \left ( \log(k+1)^0- \log(k)^0 \right )=0
\end{equation}
vanishes for all $k\geq1$, this is the convention in the Stieltjes-Hermite method, see \cite[p. 282]{BC88}. \\

Summing the right side of the sequence of functions (\ref{el10}) over the index $k\geq 1$, and using the absolute convergence of both the power series and the integral yield
\begin{eqnarray}\label{el14}
\sum_{k \geq 1}v_k(s)&=&\sum_{k \geq 1} \left (\frac{\mu(k)}{k^s}-c_{-1}\int_{k}^{k+1} \frac{1}{x^{s-\rho+1}}d x \right )  \nonumber\\
&=&\sum_{k \geq 1} \frac{\mu(k)}{k^s}-c_{-1}\int_{1}^{\infty} \frac{1}{x^{s-\rho+1}}d x \\
&=& \frac{1}{\zeta(s)}-\frac{c_{-1}}{s-\rho}  \nonumber
\end{eqnarray}
for any complex number $s = \sigma + it \in \mathbb{C}$ such that $\Re e(s) = \sigma > 1$. \\

Summing the right side of the sequence of functions (\ref{el12}) over the index $k\geq 1$, and using the absolute convergence of the power series for $\Re e(s) > 1$, it follows that the double sum can be summed in any order:
\begin{eqnarray}\label{el15}
\sum_{k \geq 1}v_k(s)&=&\sum_{k \geq 1}  \frac{\mu(k)}{k^\rho}  \sum_{n\geq 0}\frac{(-1)^n(s-\rho)^n}{n!} \log (k)^n   \nonumber\\
& & \hskip .5 in + \quad c_{-1} \sum_{k \geq 1}\sum_{n\geq 0}\frac{(-1)^{n+1}(s-\rho)^n}{(n+1)!}\left ( \log (k+1)^{n+1}-\log (k)^{n+1} \right ) \\
&=&\sum_{n\geq 0}\frac{(-1)^n(s-\rho)^n}{n!}  \sum_{k \geq 1}  \left (   \frac{\mu(k) \log(k)^{n}}{k^\rho}   - c_{-1}\frac{ \log (k+1)^{n+1}-\log (k)^{n+1}}{n+1} \right )  \nonumber \\ 
&=&\sum_{n\geq 0}\frac{(-1)^n \phi_n }{n!} (s-\rho)^n  \nonumber ,            
\end{eqnarray}
where the $n$th coefficient $\phi_n$ is written in the form
\begin{equation}\label{el16}
\phi_n= \sum_{k \geq 1}  \left (   \frac{\mu(k) \log(k)^{n}}{k^\rho}   - c_{-1}\frac{ \log (k+1)^{n+1}-\log (k)^{n+1}}{n+1} \right )  ,            
\end{equation}
for $n\geq 0$. The right side of equation (\ref{el15}) is an analytic, and absolutely convergent function for all complex numbers $s \in \mathbb{C}$.\\

Therefore, from (\ref{el14}) and (\ref{el15}), it follows that the Laurent series
\begin{equation}  \label{el17}
\frac{1}{\zeta(s)}=\frac{ c_{-1}}{s-\rho}+\sum_{n\geq 0}\frac{(-1)^{n} \phi_n }{n!} (s-\rho)^n,            
\end{equation}
where $c_{-1}=1/\zeta^{'}(\rho)$ is the residue at $s=\rho$, represents an analytic continuation of the inverse zeta function $1/\zeta(s)$ to the domain $D(\rho)=\{s:0<|s-\rho|<r\} $.\\

On the contrary, suppose that the fixed nontrivial zero $\rho=\sigma_0+it_0$ has multiplicity $m > 1$, and let
\begin{equation}\label{el18}
v_k(s)=\frac{\mu(k)}{k^s}-c_{-1}\int_{k}^{k+1} \frac{1}{x^{s-\rho+1}}d x-c_{-2}\int_{k}^{k+1} \frac{1}{x^{s-\rho+1}}d x - \cdots -c_{-m}\int_{k}^{k+1} \frac{1}{x^{s-\rho+1}}d x.
\end{equation}
Then, using the same procedure as (\ref{el10}) to (\ref{el17}), the Laurent series is
\begin{eqnarray}  \label{el19}
\frac{1}{\zeta(s)}-\frac{ c_{-m}}{(s-\rho)^m}-\cdots-\frac{ c_{-2}}{(s-\rho)^2}-\frac{ c_{-1}}{s-\rho}&=&\sum_{k\geq 1} v_{k}\\
&=&\sum_{n\geq 0}\frac{(-1)^{n} R(n) }{n!} (s-\rho)^n  \nonumber,            
\end{eqnarray}
where the $n$th term
\begin{equation}\label{el20}
R(n)= \sum_{k \geq 1}  \left (   \frac{\mu(k) \log(k)^{n}}{k^\rho}   - \left ( c_{-1}+\frac{c_{-2}}{s-\rho}+ \cdots+\frac{c_{-m}}{(s-\rho)^{m-1}} \right ) \left (\frac{ \log (k+1)^{n+1}-\log (k)^{n+1}}{n+1} \right ) \right ).            
\end{equation}
The left hand side of (\ref{el19}) is analytic (it has no poles) on a small disk $D(\rho)$ of some radius $r > 0$. But the right side has a pole at $s = \rho$ of multiplicity $m-1\geq 1$ as confirmed in (\ref{el20}). These information imply that $m = 1$.   \end{proof}                                                                                                                          		

In sypnosis, The first coefficient has the form
\begin{equation}\label{el21}
\phi_0= \sum_{k \geq 1}  \left (   \frac{\mu(k)}{k^\rho}   - c_{-1} \log (1+1/k) \right )  .            
\end{equation}
It is a constant attached to a fixed nontrivial zero $\rho=\sigma_0+it_0$. This is analogous to the definition of the Euler constant: the first Stieltjes coefficient of the series (\ref{el06}) defines the Euler constant
\begin{equation}\label{el22}
\gamma=\gamma_0(1)= \sum_{k \geq 1}  \left (   \frac{1}{k}   -\log (1+1/k) \right )  .            
\end{equation}
However, the $n$th coefficient $\phi_n=\phi_n(\rho)$ carries far more information than the Stieltjes $n$th coefficient $\gamma_n=\gamma_n(1)$. For example, it includes important information on the Mertens sum $\sum_{n\leq x}\mu(n)$. \\

The sequence of coefficients $\{ \phi_n(\rho):n\geq 0\}$ is a complicated sequence of functions of $n \geq 0$ involving the partial sums of the $k$th derivatives of the inverse zeta 
function
\begin{equation}
\frac{d^{(k)}}{ds} \zeta^{-1}(s)=(-1)^k\sum_{n\geq 1}\frac{\mu(n) \log ^k n}{n^{s}}  ,
\end{equation}
which is a topic of current research, see \cite{BS02}, \cite{DF02}, and the literature. This 
sequence is somewhat similar to the sequence of Stieltjes constants $\{\gamma_n(1):n\geq 1\}$, see \cite[Eq. 25.2.4]{DLMF}. The properties of the sequence of Stieltjes constants $\gamma_n$ are studied in \cite{AJ11}, \cite{CM09}, 
\cite{SS13}, \cite{KC11}, \cite{KJ92}, and \cite{KR03}. \\

\vskip .75 in

In these studies it is shown that the Stieltjes constants $\gamma_n=\gamma_n(1)$ are unbounded functions of $n \geq 1$, for example, each one satisfies the inequality $|\gamma_n|\leq (e \cdot n!)/(2^n \sqrt{n})$ is given in \cite{CM09}. A survey of these studies, and an improved bound and information on the sign changes are undertaken in \cite{AJ11}. The generalized Stieltjes constants to $L$-functions appears in \cite{SS13}.

\section{A Second Proof} \label{s666}
The zeta function has a myriad of different representations. Among these is the Jensen integral representation of the zeta function 
\begin{equation} \label{eq666.23}
\zeta(s)=\frac{\pi}{2(s-1)}\int_{-\infty}^{\infty}\frac{\left( 1/2+it\right )^{1-s}}{\cosh(\pi t)^2} dt,
\end{equation}
where $s \in \C$, see \cite[Theorem 1]{JB18}, facilitates another way of showing the simplicity of the zeros.
						
 \begin{thm} \label{thm666.02} Each nontrivial zero $\rho = \sigma + it$ of the zeta function $\zeta(s)$ is a simple zero. In particular, for each complex number $s \in \C$, the Taylor series at a fixed nontrivial zero $\rho_0= \sigma_0 + it_0$ is the power series
\begin{equation}\label{eq666.48}
\zeta(s)=\sum_{n \geq 1}\frac{\zeta^{(n)}(\rho_0)}{n!}\left ( z-\sigma_0- it_0\right )^n.
\end{equation}
\end{thm}
		
\begin{proof} Let $\rho = \sigma + it$ be a nontrivial zero. The change of variable $s=1+\sigma - it+z$ in (\ref{eq666.23}) leads to
\begin{equation}\label{eq666.33}
\zeta(1-\sigma- it+z)=\frac{\pi}{2(z-\sigma- it)}\int_{-\infty}^{\infty}\frac{\left( 1/2+it\right )^{-\sigma+ it-z}}{\cosh(\pi t)^2} dt,
\end{equation}
where $z \in \C$. It has a simple pole at $z = \sigma + it$. The integral 

\begin{equation}\label{eq666.36}
\left |\int_{\infty}^{\infty}\frac{\left( 1/2+it\right )^{-\sigma+ it-z}}{\cosh(\pi t)^2} dt \right | \ll \int_{0}^{\infty}\frac{1}{\cosh(\pi t)^2} dt <\infty
\end{equation}
is absolutely bounded for any complex number $z \in \C$. Evaluation at a nontrivial zero $z= \sigma + it$, implies that the simple pole on the left side 
\begin{equation}\label{eq666.38}
\zeta(1)=\frac{\pi}{2(z-\sigma- it)}\int_{-\infty}^{\infty}\frac{\left( 1/2+it\right )^{-\sigma+ it-z}}{\cosh(\pi t)^2} dt
\end{equation}
is matched with a simple pole on the right side. Accordingly, the Taylor series at a fixed nontrivial zero $\rho_0= \sigma_0 + it_0$ is
\begin{equation}\label{eq666.48}
\zeta(s)=\sum_{n \geq 1}\frac{\zeta^{(n)}(\rho_0)}{n!}\left ( s-\sigma_0- it_0\right )^n
\end{equation}
is well defined for all complex numbers. Here the $n$th derivative
\begin{equation} \label{eq666.23}
\zeta^{(n)}(s)=\frac{d^{n}}{ds^n} \left (\frac{\pi}{2(s-1)}\int_{-\infty}^{\infty}\frac{\left( 1/2+it\right )^{1-s}}{\cosh(\pi t)^2} dt \right ),
\end{equation}
is absolutely continous for all complex numbers $s\ne 1$. The evaluation at $s= \sigma_0 + it_0$ yields the $n$th Taylor coefficient $\zeta^{(n)}(\rho_0)/n!$ for every $n \geq 0$. This proves the claim.
\end{proof}

\section{Generalization Of The Stieltjes-Hermite Method} \label{s333}
The Stieltjes-Hermite method of the previous Section should be extendable to the Laurent series some other zeta functions and theirs inverses, 
and $L$-functions and theirs inverses. For example, a more general zeta function and its inverse are,
\begin{equation} \label{el29}
\zeta_{F}(s)= \sum_{\mathfrak{a} \subset \mathbb{Z}_F}\frac{1}{N(\mathfrak{a})^s}     \qquad \text{     and      } \qquad   \frac{1}{\zeta_F(s)}= \sum_{\mathfrak{a} \subset \mathbb{Z}_F}\frac{\mu_F(\mathfrak{a})}{N(\mathfrak{a})^s},
\end{equation}
where $F$ is a global number fields, and $\mathbb{Z}_F$ is its ring of integers, see \cite[p. 315]{NW04} for a generalization of the Mobius function. 
And an $L$-function and its inverse
\begin{equation}
L(s,f)= \sum_{n \geq 1}\frac{\lambda_f(n)}{n^s}   \qquad   \text{     and      }  \qquad \frac{1}{L(s,f)}= \sum_{n \geq 1}\frac{\mu_f(n)}{n^s},
\end{equation}
where $\sum_{d|n}\mu_f(d)\lambda(n/d)=0$ for $n>1$, respectively.\\

The generalization of the Stieltjes-Hermite method to the Laurent series of these functions is clearly and heavily dependent on many parameters such as the 
\begin{enumerate} [font=\normalfont, label=(\roman*)]
\item The coefficients local, global and functions fields $F$ in characteristics char $F=0$ and char $F>0$.
\item The existence of Euler Products.
\item The existence of zeros or poles at $s=1/2,1$, etc. 
\item The Galois groups $Gal(E/F)$, and theirs irreducible representations. 
\item The regulator $R_F$, the class number $h_F(d) \geq 1$ of $F$, and and the discriminant $d_{F}$.
\end{enumerate}

These parameters can vary the multiplicities of the zeros and poles of the Laurent series for these functions. Hence, these parameters can turn 
the analysis into a very complicated subject, both algebraically and analytically. It appears that the most important parameters are the 
coefficients fields $F$, the Euler products, the class number, and the Galois group $Gal(E/F)$.\\

\subsection{Abelian Zeta and $L$-Functions}
The zeta functions associated with global fields with Abelian Galois groups $Gal(E/F)$ have representations as products of linearly independent $L$-functions;
\begin{equation}
\zeta_{F}(s)= \prod_{\chi} L(s,\chi),
\end{equation}
confer \cite[p. 414]{NW04}. All the numerical evidence, and theoretical information available in the vast literature demonstrate that abelian zeta 
functions over global fields of characteristic $char(F)=0$ have simple zeros, see \cite{LMFDB}, \cite{BA12}, \cite{MN13}.  \\

\textbf{Example 3.1.} The zeta function over the rational numbers is an ideal case: it has trivial coefficients, Galois group $Gal(F/\mathbb{Q})={1}$, class number 1, and a perfect Euler product, and so on. This is also demonstrated by the trivial residue $r_{-1}=1$ of the Laurent series (\ref{el06}) at $s=1$. In contrast, the Dedekind zeta  
\begin{equation} \label{el30}
\zeta_F(s)=\frac{ \beta_{-1}}{s-1}+\sum_{n\geq 0} \frac{(-1)^n \beta_n}{n!}(s-1)^n,
\end{equation}
for a field extension $F$ of $\mathbb{Q}$ has a very complex residue 
\begin{equation}
\beta_{-1}=\frac{2^{r_1}(2\pi)^{r_2}h_{F}R_{F}}{w_{F}\sqrt{|d_{F}|}},
\end{equation}
see \cite[p. 37]{CN09}, and similar references.\ This can have significant effect on the coefficients of the inverse zeta function (\ref{el29}), and the Laurent series
\begin{equation} \label{el32}
\frac{1}{\zeta_F(s)}=\frac{ c_{-m}(F)}{(s-\alpha)^m}+ \cdots+\frac{ c_{-1}(F)}{s-\alpha}+\sum_{n\geq 0} \frac{(-1)^n c_n(F)}{n!}(s-\alpha)^n.
\end{equation}

\subsection{Nonabelian Zeta and $L$-Functions }
Let $\rho: Gal(E/F) \longrightarrow GL(n)$ be an irreducible representation of the Galois group $Gal(E/F)$ of degree $deg(\rho)=d$. The Dedekind 
zeta functions associated with fields with nonabelian Galois groups $Gal(E/F)$ have representations as products of linearly independent Artin 
$L$-functions;
\begin{equation}
\zeta_{F}(s)= \prod_{\rho} L(s,\rho)^d.
\end{equation}
The occurrence of one or more nonlinear representation of degree $deg(\rho)=d>1$  forces the zeta function to have zeros of multiplicity $d>1$. \\

\textbf{Example 3.2.} An example of a field extension $F=\mathbb{Q}(\theta)$, where $\theta$ is a root of $f(x)=x^4-2x^2+2$, which has the Galois group $Gal(F/\mathbb{Q})=D_4$, and the irreducible representation  $\rho: Gal(\overline{\mathbb{Q}}/\mathbb{Q}) \longrightarrow GL(n)$ of degree $deg(\rho)=2$, is listed in \cite{LMFDB}. The corresponding Dedekind zeta function has the representation\\
\begin{equation}
\zeta_{F}(s)= L(s,\rho_0)\prod_{\rho} L(s,\rho)^2.
\end{equation}
In this case, the Laurent series of the inverse zeta function should have the form\\ 
\begin{equation}
\frac{1}{\zeta_F(s)}=\frac{ c_{-2}}{(s-\alpha)^2}+\frac{ c_{-1}}{s-\alpha}+\sum_{n\geq 0} \frac{(-1)^n c_n}{n!}(s-\alpha)^n.
\end{equation} 
Here, the basic Stieltjes-Hermite method can fail for infinitely many nontrivial zeros $\alpha$ of $\zeta_F(s)$ of multiplicity 2. Thus, the basic Stieltjes-Hermite method  has to be modified to handle these cases. However, it is not clear how the algebraic properties of the coefficients of the Laurent series, see (\ref{el32}), obstruct the analytic properties whenever the Galois group is nonabelian.

\textbf{Example 3.3.} An example of an $L$-function of degree 2 Artin $L$-function for a $S_3$ extension $F$ of the rational numbers $\mathbb{Q}$ with square local factors in its Euler product 
\begin{equation}
L(s,\rho)= \prod_{p\geq 2} L_{p}(s, \rho)
\end{equation} was computed in \cite[Lemma 4.8]{CF13}. Introductions to the calculations of zeta and $L$-functions are given in \cite[p. 167]{SJ00}, and \cite{BA05}, and similar references.

\currfilename.\\

\begin{thebibliography}{99}
\bibitem{AJ11} Adell, J. A. \textit{\color{green}Asymptotic estimates for Stieltjes constants: a probabilistic approach}. Proc. R. Soc. Lond. Ser. A Math. Phys. Eng. Sci. 467 (2011), no. 2128, 954-963.

\bibitem{AP85} Apostol, T. M. \textit{\color{green}Formulas for higher derivatives of the Riemann zeta function}. Math. Comp. 44 (1985), no. 169, 223-232.

\bibitem{BA12} Andrew R. Booker. \textit{\color{green}Simple zeros of degree 2 L-functions}, arXiv:1211.6838.

\bibitem{BA05} Andrew R. Booker. \textit{\color{green}Artin's conjecture, Turing's method and the Riemann hypothesis}, arXiv:math/0507502.

\bibitem{BH13} Bui, H. M.; Heath-Brown, D. R. \textit{\color{green}On simple zeros of the Riemann zeta-function}, arXiv:1302.5018.

\bibitem{BC88} Bohman, Jan; Frieberg; Carl-Erik. \textit{\color{green}The Stieltjes function definition and properties.} Math. Comp. 51 (1988), no. 183, 281-289.

\bibitem{BS02} Binder, Thomas; Pauli, Sebastian; Saidak Filip. \textit{\color{green}New zero free regions for the derivatives of the Riemann zeta function,} arXiv:1002.0362.

\bibitem{CM09} Coffey, Mark W.\textit{\color{green} Series representations for the Stieltjes constants}, arXiv:0905.1111.

\bibitem{CG98} Conrey, J. B.; Ghosh, A.; Gonek, S. M. \textit{\color{green}Simple zeros of the Riemann zeta-function.} Proc. London Math. Soc. (3) 76 (1998), no. 3, 497-522.

\bibitem{CO07} Cohen, Henri.  \textit{\color{green}Number theory. Vol. II. Analytic and modern tools.} Graduate Texts in Mathematics, 240. Springer, New York, 2007.

\bibitem{CN09}  Childress, Nancy. \textit{\color{green}Class field theory}. Universitext. Springer, New York, 2009.

\bibitem{CF13} Calegari, Frank. \textit{\color{green}The Artin conjecture for some $S_5$-extensions}. Math. Ann.  356  (2013),  no. 1, 191-207. 

\bibitem{DF02} Duenez, Eduardo; Farmer, David W.; Froehlich, Sara; Hughes, Chris; Mezzadri, Francesco; Phan, Toan. \textit{\color{green}Roots of the derivative of the Riemann zeta function and of characteristic polynomials}, arXiv:1002.0372.

\bibitem{FW04} David W. Farmer, Kevin Wilson. \textit{\color{green}Converse theorems assuming a partial Euler product}, arXiv:math/0408221.  

\bibitem{GO00} Goldston, D. A.; Gonek, S. M.; Ozulik, A. E.; Snyder, C. \textit{\color{green}On the pair correlation of zeros of the Riemann zeta-function}. Proc. London Math. Soc. (3) 80 (2000), no. 1, 31-49.

\bibitem{IA12} Aleksandar Ivic. \textit{\color{green}The Theory of Hardy's Z-Function}, Series: Cambridge Tracts in Mathematics (No. 196), 2012.

\bibitem{IA03} Ivic, Aleksandar. \textit{\color{green}The Riemann zeta-function. Theory and applications}. Wiley, New York; Dover Publications, Inc., Mineola, NY, 2003.

\bibitem{DLMF} DLMF. \textit{\color{green}Digital Library Mathematical Functions}, http://dlmf.nist.gov.

\bibitem{HG11} Hiary, Ghaith Ayesh. \textit{\color{green}Fast methods to compute the Riemann zeta function.} Ann. of Math. (2)  174  (2011),  no. 2, 891-946.

\bibitem{HG16} Hiary, Ghaith A. \textit{\color{green}An alternative to Riemann-Siegel type formulas}. Math. Comp.  85  (2016),  no. 298, 1017-1032. 

\bibitem{JB18} Fredrik Johansson, Iaroslav Blagouchine. \textit{\color{green}Computing Stieltjes constants using complex integration}, arXiv:1804.01679.

\bibitem{KA93} Karatsuba, Anatolij A. \textit{\color{green}Basic analytic number theory}. Translated from the second (1983) edition. Springer-Verlag, Berlin, 1993.  

\bibitem{KC11} Knessl, Charles; Coffey, Mark W. \textit{\color{green}An effective asymptotic formula for the Stieltjes constants}. Math. Comp. 80 (2011), no. 273, 379-386.

\bibitem{KJ92} Keiper, J. B. \textit{\color{green}Power series expansions of Riemann's $\xi$ function}. Math. Comp. 58 (1992), no. 198, 765-773.

\bibitem{KR03} Kreminski, Rick. \textit{\color{green}Newton-Cotes integration for approximating Stieltjes (generalized Euler) constants}. Math. Comp. 72 (2003), no. 243, 1379-1397.

\bibitem{LM74} Levinson, Norman; Montgomery, Hugh L. \textit{\color{green}Zeros of the derivatives of the Riemann zetafunction.} Acta Math. 133 (1974), 49-65.

\bibitem{LMFDB} LMFDB. \textit{\color{green}The L-functions and modular forms database}, www.lmfdb.org.

\bibitem{MH72} Montgomery, H. L. \textit{\color{green}The pair correlation of zeros of the zeta function.} Analytic number theory (Proc. Sympos. Pure Math., Vol. XXIV, St. Louis, Mo., 1972), pp. 181-193. A. Math. Soc., Providence, R.I., 1973.

\bibitem{MN06} Milinovich, Micah B.; Ng, Nathan. \textit{\color{green}Lower bounds for moments of zeta prime rho}, arXiv:0706.2321.

\bibitem{MN13} Milinovich, Micah B.; Ng, Nathan. \textit{\color{green}Simple zeros of modular L-functions,} arXiv:1306.0854. 

\bibitem{MV07} Montgomery, Hugh L.; Vaughan, Robert C. \textit{\color{green}Multiplicative number theory. I. Classical theory}. Cambridge University Press, Cambridge, 2007.

\bibitem{NW04} Narkiewicz, Wladyslaw. \textit{\color{green}Elementary and analytic theory of algebraic numbers}. Third edition. Springer Monographs in Mathematics. Springer-Verlag, Berlin, 2004.

\bibitem{NW00} Narkiewicz, W. \textit{\color{green}The development of prime number theory. From Euclid to Hardy and Littlewood.} Springer Monographs in Mathematics. Springer-Verlag, Berlin, 2000.

\bibitem{OT85} Odlyzko, A. M.; te Riele, H. J. J. \textit{\color{green}Disproof of the Mertens conjecture.} J. Reine Angew. Math.  357  (1985), 138-160.

\bibitem{RP96} Ribenboim, Paulo. \textit{\color{green}The new book of prime number records,} Berlin, New York: Springer-Verlag, 1996.

\bibitem{SJ02} Steuding, J. \textit{\color{green}On simple zeros of the Riemann zeta-function in short intervals on the critical line}. Acta Math. Hungar. 96 (2002), no. 4, 259-308. 

\bibitem{SJ00} Steuding, J. \textit{\color{green}Introduction To $L$-Functions}, Preprint 2000, Available on www. 

\bibitem{SN10} Snaith, N. C. \textit{\color{green}Riemann zeros and random matrix theory.} Milan J. Math. 78 (2010), no. 1, 135-152.

\bibitem{SS13} Saad Eddin, Sumaia. \textit{\color{green}Explicit upper bounds for the Stieltjes constants.} J. Number Theory 133 (2013), no. 3, 1027-1044.

\bibitem{TE97} Tollis, Emmanuel. \textit{\color{green}Zeros of Dedekind zeta functions in the critical strip.} Math. Comp.  66  (1997),  no. 219, 1295-1321.

\bibitem{TE86} Titchmarsh, E. C. \textit{\color{green}The theory of the Riemann zeta-function.} Second edition. Edited and with a preface by D. R. Heath-Brown. The Clarendon Press, Oxford University Press, New York, 1986. 

\bibitem{YC00} Yildirim, Cem Yaliin. \textit{\color{green}Zeros of $\zeta^{''}(s) , \zeta^{'''}(s), s<1/2$.} Turkish J. Math. 24 (2000), no. 1, 89-108.

\bibitem{ZD75} Zagier, Don A. \textit{\color{green}Kronecker limit formula for real quadratic fields}. Math. Ann. 213 (1975), 153-184.
 
\end{thebibliography}
\end{document}